\documentclass[12pt]{article}
\usepackage{amsmath,amssymb,amscd,latexsym,mathrsfs}
\usepackage{color}

\usepackage[all]{xy}

\newcommand{\Lim}{\underleftarrow{\lim}}
\newcommand{\coLim}{\underrightarrow{\lim}}
\newcommand{\dom}{{\rm dom\,}}
\newcommand{\cod}{{\rm cod}}

\newcommand{\Cat}{{\rm Cat}}
\newcommand{\CAT}{{\rm CAT}}
\newcommand{\Set}{{\rm Set}}
\newcommand{\sSet}{{\rm sSet}}

\newcommand{\Ab}{{\rm Ab}}
\newcommand{\Ob}{{\rm Ob}}
\newcommand{\Mor}{{\rm Mor}}

\newcommand{\Imm}{{\rm Im\,}}
\newcommand{\Ker}{{\rm Ker\,}}

\newcommand{\diag}{{\rm diag\,}}

\newcommand{\FF}{{\mathbb F}}
\newcommand{\DD}{{\mathbb D}}
\newcommand{\DF}{{\mathbb{FD}}}

\newcommand{\NN}{{\,\mathbb N}}

\newcommand{\mC}{{\mathscr C}}
\newcommand{\mD}{{\mathscr D}}

\newcommand{\mF}{{\cal F}}
\newcommand{\mG}{{\cal G}}
\newcommand{\mH}{{\cal H}}

\newcommand{\ZZ}{{\,\mathbb Z}}

\newcommand{\mA}{{\mathscr A}}

\newcommand{\cA}{{\mathcal A}}

\newcommand{\mE}{{\mathcal E}}

\newcommand{\fF}{{\mathfrak F}}
\newcommand{\fU}{{\mathfrak U}}

\newtheorem{theorem}{\bf Theorem}[section]
\newtheorem{lemma}[theorem]{\bf Lemma}
\newtheorem{proposition}[theorem]{\bf Proposition}
\newtheorem{corollary}[theorem]{\bf Corollary}
\newtheorem{definition}{\sc Definition}[section]
\newtheorem{example}[definition]{\sc Example}
\newtheorem{remark}[definition]{\sc Remark}

\def\leq{\leqslant}
\def\geq{\geqslant}

\begin{document}

\begin{center}
{\large \bf
Criteria for preserving the category cohomology for the inverse image
}
\\
\medskip
Ahmet A. Husainov\\
\end{center}

\begin{abstract}
The article investigates the question of under what conditions a functor between small categories preserves cohomology groups when passing to the inverse image.
For example, it is known that the left adjoint functor preserves the category cohomology with local coefficients or the category cohomology constructed as derived of the limit functor.
We give counterexamples showing that the left adjoint functor may not preserve the Baues-Wirsching, Hochschild-Mitchell, and Thomason cohomology.
To solve the arising problems, we propose necessary and sufficient conditions for the invariance of these types of cohomology under the transition to the inverse image of a functor between small categories. Moreover, we generalize these cohomology of small categories and find similar criteria for the obtained generalization.
\end{abstract}

2010 Mathematics Subject Classification 55U10, 18G10, 18G30, 18G35


\tableofcontents

\section{Introduction}

The cohomology of small categories appeared as a means to measure the deviation from the exactness property of the limit functor.
It was used to calculate the cohomological dimension, to classify extension groups and other constructions, to study spaces with local coefficient systems.
After the appearance of the paper \cite{bau1985}, the scope of applications of cohomology of small categories has expanded.
The Baues-Wirsching cohomology has been applied to the study of algebraic theories \cite{jib2006}.
A spectral sequence was constructed for the Baues-Wirsching cohomology of the Grothendieck construction \cite{pir2006}, and it was developed in \cite{gal2012}.
The paper \cite{gal2021} was devoted to the homology theory of Gabriel-Zisman \cite{gab1967}. Applications of the category cohomology theory to the theory of extensions of small categories are developed in \cite{yal2023}. The theory of cohomology of categories of subgroups in a discrete group is developed in \cite{yal2022}.

Conditions for the invariance of cohomology in passing to the inverse image have always helped to solve problems in the theory of (co)homology of small categories.
For example, to prove that the cohomological dimension of the directed cofinality set
 $\aleph_n$ is equal to $n+1$ \cite[Theorem A]{mit1973}, for Mitchell, it was enough to prove the cohomology conservation theorem when passing to the inverse image along the cofinal functor \cite[Theorem B]{mit1973}, and to use the formula for the cohomological dimension of a totally ordered set \cite[Corollary 36.9]{mit1972}.
We remark, that Mitchell's investigations on the cohomological dimension 
of totally ordered sets \cite{mit1972} are continued in \cite{ber2021}.

For the author of this paper, the study of invariance conditions for coho\-mo\-lo\-gy of categories with coefficients in natural systems \cite{bau1985} when passing to the inverse image was necessary to characterize small categories with cancellations whose Hochschild-Mitchell dimension is at most 1 \cite{X1997}.

It is known that the cohomology of small categories with coefficients in diagrams of abelian groups is invariant under the inverse image functor along a functor admitting a right adjoint.
In \cite[Theorem 5.4 and Theorem 5.10]{mur2006}, sufficient conditions were obtained for the preservation of the Baues-Wirsching cohomology for some natural systems and adjoint functors between small categories.

 The Baues-Wirsching cohomology is isomorphic to the cohomology of the category of factorizations. In this connection, there have been reports in the press that the Baues-Wirsching cohomology and the Thomason cohomology have the property of invariance with respect to the left adjoint functor.
We give counterexamples and describe necessary and sufficient conditions for functors along which the inverse image functor preserves the Baues-Wirsching, Thomason, and Hochschild-Mitchel cohomology values.

At the end of the paper, we introduce compatible cohomology for small categories generalizing these cohomology and find similar criteria for compatible cohomology.

\section{Notation and Preliminaries}

Let us write out the initial definitions and notation. The rest will be given in the course of the presentation.

\begin{itemize}
\item $\Set$ - category of sets and mappings.
\item For any (locally small) category $\mA$, we denote the set of morphisms $a\to b$ between $a, b\in \Ob\mA$ by $\mA(a,b)$, and the morphism functor - via $\mA(-, =): \mA^{op}\times \mA\to \Set$. If $\mA$ is an Abelian category, then by 
$Hom_{\mA}(a,b)$ or $Hom(a,b)$ we denote the Abelian morphism group.
\item
$\mA^{\mC}$ is the category of functors from the small category $\mC$ to an arbitrary category $\mA$.
Functors $F$ from a small category $\mC$ to an arbitrary $\mA$ are called object diagrams of the category $\mA$ over $\mC$ and can be denoted as the family $\{F(c)\}_{c\in \mC}$.
\item $\Delta_{\mC}A$ is a functor $\mC\to \mA$ that takes the constant values 
$A\in \Ob\mA$ on objects and the values $1_A$ on morphisms.
\item $\NN$ - set of non-negative integers.
\item $\Cat$ - category of small categories and functors. We call a category small 
if the set of its morphisms belongs to some universe
$\fU$ (see \cite[\S I.1]{mac2004}). It is assumed that the universe $\fU$ is larger than the universe of finite sets.
\item $\CAT$ - category of categories whose morphism classes are subsets of the universe 
$\fU$.
\item $h_c= h^{\mC}_c$, for the small category $\mC$ and its object $c\in \Ob\mC$, denotes the morphism functor, $h_c(-): \mC^{op} \to \Set$.
\item $0$ - a group consisting of one element.
\item $\Delta$ - the category of finite linearly ordered sets $[n]= \{0, 1, \cdots, n\}$, 
$n\geq 0$, and non-decreasing mappings. The category $\Delta$ is generated by 
morphisms of the following form:
\begin{enumerate}
\item $\partial^i_n: [n-1]\to [n]$ (for $0\leq i\leq n$) - increasing mapping whose image does not contain $i$,
\item
$\sigma^i_n: [n+1]\to [n]$ (for $0\leq i\leq n$) is a nondecreasing surjection that takes the value $i$ twice.
\end{enumerate}
\end{itemize}

Let $\Phi: \mC\to \mD$ be a functor between categories.
For each $d\in \Ob\mD$ the left fibre \cite[Appendix 2, \S3.5]{gab1967} or the comma category 
$\Phi$ over $d$ \cite{mac2004} is the category $\Phi\downarrow d$, whose objects are the pairs 
$(c\in \Ob\mC, \beta\in \mD(S(c), d))$, and the morphisms $(c, \beta)\to (c', \beta')$ are 
given by $\alpha\in \mC(c,c')$ satisfying $\beta'\circ S(\alpha)= \beta$. A morphism in 
$\Phi\downarrow d$ is defined by a triple consisting of two objects and a morphism between them.

For an arbitrary category $\mA$, the inverse image functor $\Phi^*: \mA^{\mD}\to \mA^{\mC}$ is defined, which assigns to each diagram $F\in \mA^{\mD}$ the
composition 
$F\Phi= F\circ \Phi\in \mA^{\mC}$, and to each natural transformation $\eta: F\to F'$ the natural transformation $\eta\Phi: F\Phi \to F'\Phi$ defined by the formula $(\eta\Phi)_{c}= \eta_{\Phi(c)}$, for all $c\in \Ob\mC$.
If $\mA$ is a cocomplete category, then the functor $\Phi^*$ has a left adjoint functor $Lan^{\Phi}: \mA^{\mC}\to \mA^{\mD}$, which is called the left Kan extension.

For an arbitrary diagram $F\in \mA^{\mC}$, the diagram $Lan^{\Phi}F$ can be considered as taking values on the objects $d\in \Ob\mD$ equal to $Lan^{\Phi}F (d)= \coLim^{\Phi\downarrow d}FQ_d$, where
$Q_d: \Phi\downarrow d\to \mC$ assigns to each $S(c)\to d$ an object $c\in \mC$, with an obvious extension to morphisms
\cite[\S10.3 (10)]{mac2004}.

A simplicial set is a functor $X: \Delta^{op}\to \Set$. The values of this functor on morphisms generating the category $\Delta$ are denoted by $d^n_i= X(\partial^i_n)$, $s^n_i= X(\sigma^i_n)$, for all $0\leq i\leq n$.
A simplicial mapping $X\to Y$ between simplicial sets is a natural transformation.
The category of simplicial sets is denoted by $\Set^{\Delta^{op}}$ or $\sSet$.

The homology groups $H_n(X)$ of a simplicial set are defined by the formula 
$H_n(X)= \Ker d_n/ \Imm d_{n+1}$ as the homology groups of a chain complex
$$
0 \leftarrow C_0(X) \xleftarrow{d_1} C_1(X) \xleftarrow{d_2} C_2(X)
 \xleftarrow{d_3}\cdots
$$
consisting of free abelian groups $C_n(X)= \ZZ X_n$ with bases $X_n$ and homo\-mor\-phisms $d_0=0$ and $d_n: C_n(X)\to C_{n-1}(X)$ for all $ n\geq 1$ defined on the basis elements $x\in X_n$ as $d_n(x) = \sum\limits^n_{i=0}(-1)^i d^n_i(x)$.
By Eilenberg's theorem \cite[Appendix 2, Theorem 1.1]{gab1967}, the homology groups 
of a simplicial set are isomorphic to the homology groups of its geometric realization.

Let $\mC$ be a small category. For $n\in \NN$, a path of length $n$ in $\mC$ is a finite sequence of morphisms $c_0\stackrel{\alpha_1}\to c_1\to \cdots \to c_{n-1}\stackrel{\alpha_n}\to c_n$ in the category $\mC$. For an arbitrary path of length $n\geq 1$ and an integer $i$ from the interval $0\leq i\leq n$ denote by $c_0 \stackrel{\alpha_1}\to c_1\to \cdots \to \widehat{c_i} \to \cdots \to c_{n-1} \stackrel{\alpha_n}\to c_n$ path equals
$$
\left\{
\begin{array}{ll}
c_1 \stackrel{\alpha_2}\to c_2\to \cdots \to c_{n-1}
   \stackrel{\alpha_n}\to c_n~, & \mbox{ if } i=0,\\
c_0 \stackrel{\alpha_1}\to c_1\to \cdots \to
c_{i-1} \stackrel{\alpha_{i+1}\circ\alpha_i}\longrightarrow c_{i+1} \to \cdots \to c_{n-1}
   \stackrel{\alpha_n}\to c_n~, & \mbox{ for } 1\leq i\leq n-1,\\
c_0 \stackrel{\alpha_1}\to c_1\to \cdots \to c_{n-2}
   \stackrel{\alpha_{n-1}}\to c_{n-1}~, & \mbox{ if } i=n.
\end{array}
\right.
$$
A nerve or classifying space of a small category $\mC$ is a simplicial set $B\mC: \Delta^{op}\to \Set$ equal to the restriction of the functor $\Cat(-, \mC)$ to the subcategory $\Delta\subset \Cat$.

The nerve of a category can be defined as a sequence of sets of $n$-simplices $B_n\mC= \Cat([n], \mC)$, boundary operators $d^i_n= \Cat(\partial^i_n, \mC): 
B_n\ mC\to B_{n-1}\mC$, $0\leq i\leq n$, and degeneration operators
$s^i_n= \Cat(\sigma^i_n, \mC): B_n\mC\to B_{n+1}\mC$, $0\leq i\leq n$.

Since the $n$-simplices $\sigma: [n]\to \mC$ are functors, they can be given as paths $c_0\stackrel{\alpha_1}\to c_1\to \cdots \to c_{n-1} \stackrel{\alpha_n}\to c_n$, consisting of objects $c_i= \sigma(i)$ for $0\leq i\leq n$ and morphisms $\alpha_i= \sigma(i-1<i):\sigma (i-1)\to \sigma(i)$ for $1\leq i\leq n$.

The boundary operators $d^i_n: B_n\mC\to B_{n-1}\mC$, $0\leq i\leq n$, will act as
$$
d^i_n(c_0 \stackrel{\alpha_1}\to \ldots \stackrel{\alpha_n}\to c_n)=
(c_0 \stackrel{\alpha_1}\to c_1\to \cdots \to \widehat{c_i} \to \cdots \to c_{n-1} \stackrel{\alpha_n}\to c_n),
$$
and the degeneration operators $s^i_n: B_n\mC\to B_{n+1}\mC$, $0\leq i\leq n$, -
$$
s^i_n(c_0\xrightarrow{\alpha_1} \ldots \xrightarrow{\alpha_n} c_n) = (c_0\xrightarrow{\alpha_1} \ldots \xrightarrow{\alpha_{i}} c_i \xrightarrow{id_{c_i}} c_i\xrightarrow{\alpha_{i+1}} \ldots \xrightarrow{\alpha_n} c_n).
$$

We define the integral homology groups of a small category $\mC$ as the homology groups of its nerve $H_n(B\mC)$.

\section{Derived limit functors}

This section is devoted to the cochain complex of Abelian groups whose cohomology groups are isomorphic to the values of the derived limit functors. Examples are given of calculating the cohomology groups of small categories with coefficients in the system of Abelian groups. Homology of small categories and left derived functors of the left Kan extension are considered.

\subsection{Definition of category cohomology}

Let $\mC$ be a small category.
The limit of the diagram $\mG: \mC\to \Ab$ is denoted by $\Lim_{\mC}\mG$.
  The category of diagrams of Abelian groups $\Ab^{\mC}$ is Abelian and has a sufficient number of injective and projective objects. Hence, for all $n\geq 0$, there are the derived limit functors 
$\Lim^n_{\mC}: \Ab^{\mC}\to \Ab$.

\begin{definition}
The cohomology of the small category $\mC$ with coefficients in the abelian group diagram $\mG$ is called the cohomology $H^n(\mC, \mG)$, $n\geq 0$, of the cochain complex $C^*(\mC, \mG)$ consisting of abelian groups
$$
0 \to C^0(\mC, \mG) \xrightarrow{d^0} C^1(\mC, \mG) \xrightarrow{d^1}
C^2(\mC, \mG) \xrightarrow{d^2} \cdots ,
$$
given as
$$
  C^n(\mC, \mG)= 
 \prod_{c_0\xrightarrow{\alpha_1} c_1\to  \ldots
 \to c_{n-1}\xrightarrow{\alpha_n}c_n}\mG(c_n)
$$
and differentials $d^n= \sum^{n+1}_{i=0} (-1)^i d^n_i$, where
$d^n_i: C^n(\mC, \mG)\to C^{n+1}(\mC, \mG)$ are coboundary operators given 
by the formula
\begin{multline*}
d^n_i\varphi(c_0\xrightarrow{\alpha_1} c_1\to  \ldots
 \to c_{n}\xrightarrow{\alpha_{n+1}}c_{n+1})=\\
 \begin{cases}
\varphi(c_0\xrightarrow{\alpha_1} c_1\to \ldots \to \widehat{c_i}\to \ldots \to c_{n}
 \xrightarrow{\alpha_{n+1}} c_{n+1}), & 0\leq i \leq n,\\
 \mG(\alpha_{n+1})
  (\varphi(c_0\xrightarrow{\alpha_1} c_1\to  \ldots \xrightarrow{\alpha_{n}} c_{n})), 
  & i = n+1. 
 \end{cases}
\end{multline*}
\end{definition}

\begin{proposition}
For any small category $\mC$, a diagram $\mG: \mC\to \Ab$ and $n\geq 0$, the groups $H^n(\mC, \mG)$ are isomorphic to the values of the derived limit functors $\Lim ^n_{\mC}\mG$.
\end{proposition}
{\sc Proof} is given in \cite[Lemma 2]{oli1994}.
This assertion also follows from the more general but dual assertion proved in \cite[Appendix 2, Proposition 3.3]{gab1967} for object diagrams of the abelian category with exact coproducts.

\subsection{Examples of category cohomology calculation}

The category $\mC$ has no retractions if the composition of any morphisms not equal to the identity morphism is not equal to the identity morphism.
In this case, there exists a subcomplex $C^*_+(\mC,\mG)$ of the complex $C^*(\mC, \mG)$ which, for $n>0$, consists of products over sequences of morphisms not equal to the identity
$$
C^n_+(\mC, \mG)= \prod\limits_{c_0\xrightarrow[\not=]{} c_1 \xrightarrow[\not=]{} c_2
\xrightarrow[\not=]{} \cdots \xrightarrow[\not=]{} c_n} \mG(c_n)
$$
and equal to $C^0(\mC, \mG)$ for $n=0$.
The cohomology of this subcomplex is isomorphic to $H^n(\mC, \mG)$ \cite[Proposition 2.2]{X1997}.

\begin{example}\label{mon1}
(Idempotent monoid cohomology) Let $\mE$ be a category consisting of a unique object $1$ and two morphisms: the identity $id: 1\to 1$ and the morphism $e:1\to 1$.
And let the composition $e\circ e$ equal $e$.
Consider the functor $\cA: \mE\to \Ab$.
A category that has a single object can be interpreted as a monoid, and a diagram of abelian groups on it can be interpreted as a monoid acting on an abelian group.
Denote $A= \cA(1)$, $\varepsilon = \cA(e): A\to A$. For each $n\geq 0$, the product $C^n_+(\mE, \cA)$ consists of a single factor equal to $A$. The homomorphism $d^n$ is defined as
\begin{multline*}
d^n\varphi(1\xrightarrow{e} 1 \to \cdots \to  1\xrightarrow{e} 1)\\
= \sum\limits^n_{i=0} (-1)^i\varphi(1\xrightarrow{e} 1 \to \cdots \to  1\xrightarrow{e} 1)\\
+ (-1)^{n+1}\cA(e)(\varphi(1\xrightarrow{e} 1 \to \cdots \to  1\xrightarrow{e} 1))
\end{multline*}

The value of the function $\varphi(1\xrightarrow{e} 1 \to \cdots \to 1\xrightarrow{e} 1)$ from a set of $n$ arrows is denoted by $\varphi(e^n)$.
Get
$$
d^n\varphi(e^{n+1})=
\begin{cases}
\varphi(e^n) - \cA(e)\varphi(e^n), & \text{ if $n$ is even};\\
\cA(e)\varphi(e^n), & \text{ if $n$ is odd}.
\end{cases}
$$
Hence, the complex $C^*_+(\mE, \cA)$ consists of groups and homomorphisms
$$
0 \to A \xrightarrow{id-\varepsilon} A \xrightarrow{\varepsilon}
A \xrightarrow{id-\varepsilon} A \to \cdots
$$
The group $H^0(\mE, \cA)= \Ker d^0= \{x\in A | x=\varepsilon(x)\}$ is equal to the subgroup consisting of the fixed elements of the group $A$.
For $n>0$, the cohomology groups of $H^n(\mE,\cA)$ are equal to $0$, because for each $n\geq 0$ the element of the kernel of the homomorphism $d^{n+1}$ belongs to the image of the homomorphism $d ^n$
due to the following two implications:
$
\varepsilon(x)=0 \Leftrightarrow x=(id-\varepsilon)(x); \quad
\varepsilon(id-x)=0 \Leftrightarrow x=\varepsilon(x).
$
\end{example}

The complex $C^*_+(\mC, \mG)$ is also suitable for calculating the cohomology of partially ordered sets.

\begin{example}\label{pos1}
Consider the category $V$ consisting of three objects and two morphisms
$$
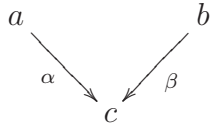
\xymatrix{
a\ar[rd]_{\alpha} && b\ar[ld]^{\beta}\\
& c
}
$$
The category $V$ also contains three identity morphisms.
Let $\cA: V \to \Ab$ be a diagram of Abelian groups taking the values $\cA(c)=A$ and $\cA(a)= \cA(b)= 0$.
To compute the groups $H^n(V, \mG)$, we construct the complex $C^*_+(V, \cA)$. There are no paths of length $\geq 2$ in the category that consist of morphisms that are not equal to the identity ones.
Hence $C^n_+(V, \cA)=0$ for $n\geq 2$. We get the complex
$$
0 \to \cA(a)\times \cA(b) \times \cA(c) \xrightarrow{d^0} \cA(c)\times \cA(c)\to 0.
$$
with differential defined as $d^0\varphi(a\to c)= \varphi(c)- \cA(\alpha)\varphi(a)$ and $d^0\varphi(b\to c)= \varphi(c)-\cA(\beta)\varphi(b)$.
Substituting $\cA(c)=A$, $\cA(a)= \cA(b)= 0$, we get $d^0$ equal to the homomorphism $\diag(\varphi)= (\varphi, \varphi) $,
and come to the complex
$$
0 \to A \xrightarrow{\diag} A\times A \to 0.
$$
Hence $H^0(V, \cA)= \Ker(\diag)= 0$, $H^1(V, \cA)= A\times{A}/\Imm(\diag)= A$ , $H^n(V, \cA)= 0$ for $n\geq 2$.
\end{example}

\subsection{Category homology groups}

\begin{definition}\label{acyclic}
A small category $\mC$ is called acyclic if $H_n(B\mC)=0$ for all $n>0$. It is called connected if $H_0(B\mC)\cong \ZZ$.
Thus, the connected category is not empty.
\end{definition}

\begin{definition}\label{defhomolcat}
The chain complex $C_*(\mC, F)$ consists of abelian groups
\begin{displaymath}
  C_n({\mC},F) = \bigoplus_{c_0 \rightarrow \cdots \rightarrow c_n} F(c_0), \quad n \geq 0,
\end{displaymath}
and homomorphisms $d_n= \sum\limits_{i=0}^{n}(-1)^i d^n_i: C_n(\mC,F) \rightarrow C_{n-1}(\mC,F)$, $n>0$,
where $d^n_i$ are defined on elements of terms
$$
(c_0 \stackrel{\alpha_1}\rightarrow c_1 \stackrel{\alpha_2}\rightarrow \cdots \stackrel{\alpha_{n}}\rightarrow c_{n}, a)
  \in \bigoplus_{c_0 \rightarrow \cdots \rightarrow c_n}F(c_0), \quad a\in F(c_0)
$$
according to the formula
$$
d^n_i(c_0 \rightarrow \cdots \rightarrow c_n, a) =
\left\{
\begin{array}{ll}
(c_0 \rightarrow \cdots \rightarrow \widehat{c_i} \rightarrow
\cdots \rightarrow c_n, a), & 1 \leq i \leq n\\
(c_1 \rightarrow \cdots \rightarrow c_n, F(c_0\stackrel{\alpha_1}\rightarrow c_1)(a) ),
& i = 0.
\end{array}
\right.
$$
Homology groups $H_n(\mC, F)$ of the small category $\mC$ with coefficients in the functor
  $F:\mC\rightarrow \Ab$ are the $n$-th homology groups of the chain complex $C_*(\mC,F)$.
\end{definition}
For example, $H_n(\mC, \Delta_{\mC}\ZZ)\cong H_n(B\mC)$, which gives grounds to call the homology of a nerve the integer homology of a category.

 \begin{proposition}
For all $n\in \NN$, the values of the left derived functors $\coLim^{\mC}_n\mG$ of the colimit functor $\coLim^{\mC}: \Ab^{\mC}\to \Ab$ are isomorphic to the homology groups $H_n(\mC, \mG)$.
  \end{proposition}
  This follows from a more general assertion about satellites of the colimit functor
  \cite[Appendix 2, Proposition 3.3]{gab1967}.
 
  Let $f: \mC\to \mD$ be a functor between small categories.
  The left Kan extension functor $Lan^f: \Ab^{\mC}\to \Ab^{\mD}$ is right exact and has left derived functors
  $Lan^{f}_n: \Ab^{\mC}\to \Ab^{\mD}$.

 According to \cite[\S10.3]{mac2004}, for any diagram $\mF: \mC\to \Ab$ the values of the left Kahn extension along $f$ to $d\in \mD$ are
  $$
   Lan^f \mF(d) = \coLim^{f\downarrow d} \mF Q_d,
  $$
where $Q_d: f\downarrow d\to \mC$ is the forgetting functor of the left fibre.
  The functor $Q_d$ maps each object $\beta: f(c)\to d$ to the object $c\in \mC$. The morphisms of the category $f\downarrow d$ are the triples of morphisms
  $(c\xrightarrow{\alpha} c', f(c)\xrightarrow{\beta} d, f(c')\xrightarrow{\beta'} d)$,
   making the triangle commutative
  $$
  \xymatrix{
  f(c)\ar[rd]_{\beta} \ar[rr]^{f(\alpha)} && f(c')\ar[ld]^{\beta'} \\
  &d
  }
  $$
  The functor $Q_d$ assigns to each such morphism a morphism $c\xrightarrow{\alpha} c'$ of the category $\mC$.
  
  \begin{proposition}\cite[Appendix 2, note 3.8]{gab1967}
For any diagram $\mF: \mC\to \Ab$ the values on $d\in \mD$ of the left derived functors 
of the left Kan extension of the diagram $\mF$ equal
  $$
  Lan^f_n \mF (d) = \coLim^{f\downarrow d}_n \mF Q_d
  $$
  \end{proposition}

\section{Category cohomology for the inverse image}

Let $f: \mC\to \mD$ be a functor between small categories.
For an arbitrary diagram of Abelian groups $\mG: \mD\to \Ab$
its inverse image is the diagram $f^*(\mG)= \mG\circ f : \mC\to \Ab$.
The term inverse image comes from the theory of sheaves.
 
This section is devoted to the spectral sequence converging to inverse image
 cohomology for a functor between small categories.
Oberst in \cite{obe1967} discussed in detail the method of constructing this spectral sequence
 for homology.
We construct it for cohomology, and mainly study the conditions for its degeneration into edge isomorphisms.

For Abelian categories $\mA$, Abelian morphism groups will be denoted by 
$Hom_{\mA}(A,B)$ or $Hom(A,B)$.

\subsection{Oberst spectral sequence}

Let $\mA$ be an abelian category with enough projective and injective objects. Then the morphism bifunctor $\mA(-,=): \mA^{op}\times \mA \to \Ab$ is balanced, in the sense that
 the functor $\mA(-,I): \mA^{op} \to \Ab$ is exact for the injective object $I$ and the functor 
$\mA(P,-): \mA \to \Ab$ is exact for the projective object $P$.

Since $\Ab^{\mC}$ is Abelian, the functor $Hom_{\Ab^{\mC}}(-,=): \left(\Ab^{\mC}\right)^{op} \times \Ab^{\mC}\to \Ab$ is balanced. Moreover, the category $\Ab$ has exact products and enough injective objects, which implies, according to \cite[III.7]{obe1967}, the existence of the following spectral sequence.
It is essential for our work, so we present the proof.

 \begin{proposition}\label{oberst1967}
For any functor $f: \mC\to \mD$ between small categories and diagrams $\mF:\mC\to \Ab$ and $\mG: \mD\to \Ab$ there exists a spectral sequence of the first quadrant
$$
E^{p,q}_2= Ext^p_{\Ab^{\mD}}(Lan^{f}_q\mF, \mG) \Rightarrow_p Ext^n(\mF, \mG\circ f)
$$
\end{proposition}
{\sc Proof.} Consider the projective resolution $\mF_*\to \mF$ of the diagram $\mF$ and the injective resolution $\mG\to \mG^*$ of the diagram $\mG$.
Let's build a double complex
$$
K^{pq}= Hom_{\Ab^{\mD}}(Lan^{f}\mF_p, \mG^q) \cong
Hom_{\Ab^{\mC}}(\mF_p, \mG^q\circ f).
$$
There are spectral sequences converging to a common limit
$I^{pq}_2= H^p H^q K^{**} \Rightarrow H^n(Tot K)$ and
$II^{pq}_2= H^q H^p K^{**} \Rightarrow H^n(Tot K)$.

  Let us calculate $II^{pq}_2$.
  
  $H^p K^{*q}= H^p(Hom(Lan^f \mF_*, \mG^q))$. The balance of $Hom$ leads to the isomorphism $H^p K^{*q} \cong Hom(Lan^f_p\mF, \mG^q)$ .
  We get
  $$
  II^{pq}_2= H^q(Hom(Lan^f_p\mF, \mG^*))= Ext^q(Lan^f_p\mF, \mG),
  $$
  according to the definition of $Ext^q$.
 
  Compute $I^{p,q}_2$.
   \begin{multline*}
  H^q K^{p *}= H^q(Hom(\mF_p, \mG^*\circ f))\\
  \cong
  Hom(\mF_p, H^q(\mG^*\circ f))=
  \begin{cases}
  Hom(\mF_p, \mG\circ f), & \text{ if $q=0$,}\\
  0, & \text{ otherwise}
  \end{cases}
  \end{multline*}
  and finally $I^{p 0}_2= H^p H^0 K^{**}= Ext^p(\mF, \mG\circ f)$, and $I^{pq}_2= 0$ if $q>0$.
  Thus, the spectral sequence $I$ degenerates, and $II$ becomes the desired one.
  \hfill$\Box$
 
 In particular, substituting $\mF= \Delta_{\mC}\ZZ$ and 
$Ext^n(\Delta_{\mD}\ZZ, \mG\circ f)= H^n(\mC, \mG\circ f)$, 
we arrive at the following statement.
 
\begin{proposition}\label{spobe}
For any functor between the small categories $f: \mC\to \mD$ and the diagram 
$\mG: \mD\to \Ab$ there exists a spectral sequence of the first quarter, of the type
$$
  E^{p,q}_2(\mG) = Ext^p(Lan^f_q\Delta_{\mC}\ZZ, \mG) \Rightarrow H^{p+q}(\mC, \mG\circ f ).
$$
\end{proposition}

\subsection{Cohomology domain of small categories}

Let $f: \mC\to \mD$ be a functor between small categories.
It gives the functor $f^*: \Ab^{\mD}\to \Ab^{\mC}$ defined on objects as $f^*(\mG)= \mG\circ f$. The functor $f^*$ assigns to each natural transformation $\xi: \mG \to \mG'$ a natural transformation $\xi*f: \mG\circ f\to \mG'\circ f$ defined on the objects $c\in \Ob\mC$ by the formula
$$
   (\xi*f)_c = \xi_{f(c)}: \mG(f(c)) \to \mG'(f(c)).
$$

Consider the functor $\Ab^{(-)}: \Cat^{op} \to \CAT$, which assigns to each small category 
$\mC$ the category $\Ab^{\mC}$ and to each functor $f: \mC \to \mD$  a functor $f^*: \Ab^{\mD}\to \Ab^{\mC}$.

We adhere to the definition of the Grothendieck construction described in \cite{tho1979}.
The Grothendieck construction $\Cat^{op}\int \Ab^{(-)}$ for the functor $\Ab^{(-)}: \Cat^{op}\to \CAT$ is a category whose objects are pairs $(\mD, \mG)$ consisting of the small category $\mD$ and the diagram $\mG: \mD\to \Ab$.
Morphisms $(f, \xi): (\mD,\mG)\to (\mC, \mF)$ in the category $\Cat^{op}\int \Ab^{(-)}$ are given as pairs, consisting of the functor $f: \mC\to \mD$ and the natural transformation $\xi: \mG\circ f \to \mF$. Composition of morphisms
  $$
  (\mE, \mH) \xrightarrow{(g, \eta)} (\mD, \mG) \xrightarrow{(f, \xi)} (\mC, \mF)
  $$
equals $(g\circ f, \xi\cdot(\eta*f) ): (\mE, \mH)\to (\mC, \mF)$. Here the dot denotes the usual composition of natural transformations.

Using the fact that $H^n(\mD, \mG)= Ext^n(\Delta_{\mD}\ZZ, \mG)$, we define homomorphisms $H^n(f, id): H^n( \mD, \mG)\to H^n(\mC, \mG\circ f)$
as compositions
\begin{multline*}
H^n(\mD, \mG)= Ext^n(\Delta_{\mD}\ZZ, \mG) \xrightarrow{Ext^n(\lambda , \mG)}
Ext^n(Lan^f\Delta_{\mC}\ZZ, \mG)=\\
  E^{n,0}_2(\mG)\xrightarrow{\gamma_n}
Ext^n(\Delta_{\mC}\ZZ, \mG\circ f) = H^n(\mC, \mG\circ f).
\end{multline*}
Here $\gamma_n$ is the edge morphism of the spectral sequence from Proposition \ref{spobe}, and $\lambda: Lan^f\Delta_{\mC}\ZZ\to
\Delta_{\mD}\ZZ$ - natural transformation,
defined as follows.
The Abelian group $Lan^f\Delta_{\mC}\ZZ(d)= H_0(f\downarrow d)$ is freely generated by the connected components of the category $f\downarrow d$ and $\lambda_d: H_0(f\downarrow d)
\to \ZZ$ maps each connected component of the category $f\downarrow d$ to the element $1\in \ZZ$.

\begin{proposition}
The cohomology of small categories given as derived limit functors on $\Ab^{\mC}$ for each small category $\mC$ can be extended to the functor $H^n(-,=): \Cat^{op}\int \Ab^{(-)}\to \Ab$ on the Grothendieck category, for all $n\geq 0$.
\end{proposition}
{\sc Proof.} For each category $\mC$ the functor $H^n(Id_{\mC}, -): \Ab^{\mC}\to \Ab$ will be equal to the functor $Ext^n(\Delta_{\mC}\ZZ,-)$. This gives the definition of $H^n(Id_{\mC}, \xi)$ for natural transformations of $\xi$ diagrams over $\mC$.
Hence, the functor 
$$
H^n(-,=): \Cat^{op}\int \Ab^{(-)}\to \Ab
$$ 
must be defined on the objects $(\mD, \mG)$ as $ H^n(\mD, \mG)$, and on the morphisms $(f, \xi): (\mD, \mG)\to (\mC, \mF)$ as compositions
$$
H^n(f, \xi)= H^n(Id_{\mC}, \xi)\circ H^n(f, id).
$$
\hfill$\Box$

Thus, the cohomology domain of small categories
is the Grothendieck construction $\Cat^{op}\int \Ab^{(-)}$.

\subsection{Verdier and Oberst criteria}

Oberst \cite{obe1967} describes a theorem proved by Verdier: If the comma categories $f\downarrow d$ are connected and acyclic for each $d\in \Ob\mD$, then the homomorphisms $H^n(f, id) $ are isomorphisms for for any diagram of Abelian groups $\mG$ and $n\geq 0$.

Later it was found that the converse is also true.

A similar statement is also true for homology.
Oberst proved \cite[Theorem 2.3]{obe1968} that for a functor between small categories $f: \mC\to \mD$ the canonical homomorphisms
$$
\coLim^{\mC}_n \mG f \to \coLim^{\mD}_n\mG
$$
are isomorphisms for any diagram $\mG: \mD\to \Ab$ if and only if for each $d\in \Ob\mD$ the category $d\downarrow f$ is connected (and therefore non-empty) and acyclic.

We will find out under what necessary and sufficient conditions these canonical homomorphisms will be isomorphisms for $n= 0, 1, 2, \cdots, N$, given $N\in \NN$.
For this purpose, consider the family of spectral sequences lying in the first quadrant (in particular, it can be a spectral functor, in the sense of \cite[\S II.2.4]{gro1957})
$$
E^{p,q}_2(\mG)\Rightarrow_p H^n(\mG),
$$
indexed by the class of objects of the Abelian category $\mA$.
For example, $\mA= \Ab^{\mD}$ and the family values belong to $\Ab$.

We are interested in the following properties that this family can have for a given natural number $N$:

\begin{itemize}
\item [($A_N$)]
For each $\mG\in \mA$, the edge morphisms $E^{n,0}_2(\mG) \to H^n(\mG)$ are  isomorphisms for all $0\leq n\leq N$ and monomorphism for $n=N+1$.
  \item[$(B_N)$]
  For each $\mG\in \mA$, $E^{0q}_2(\mG)=0$ is true for all $q$ such that $1\leq q\leq N$.
\end{itemize}
Property $(A_0)$ is always satisfied.
If $(A_N)$ is true, then we say that the spectral sequence degenerates into edge isomorphisms at the $N$ level. If $(A_N)$ is true for all $N\in \NN$, then it is called degenerate.

\begin{lemma}\label{degABN}
Let $N\geq 1$ be a natural number.
Let for each object $\mG$ of the Abelian category $\mA$ there be given a spectral sequence of Abelian groups lying in the first quadrant $E^{p,q}_2(\mG)\Rightarrow_p H^n(\mG)$ satisfying the following condition:
If $E^{0q}_2(\mG)=0$ for each $\mG$ and for all $1\leq q\leq N$, then $E^{pq}_2(\mG)=0$ for each $\mG$, for all $p\geq 1$ and for all $1\leq q\leq N$. (In this case, strings with numbers $1, 2, \ldots, N$ will consist of zero groups.)

Then $(A_N)$ is satisfied if and only if $(B_N)$ is satisfied.
\end{lemma}

\begin{remark}
The conditions of the Lemma \ref{degABN} are satisfied by the spectral sequence from Proposition \ref{oberst1967}: If $Hom(Lan^f_q\mF, \mG)=0$ for some $q$ for all $\mG\in \Ab ^{\mD}$, then $Lan^f_q\mF=0$, and hence $Ext^p(Lan^f_q\mF, \mG)=0$ for all $p\geq 0$ and $\mG\ in \Ab^{\mD}$.
\end{remark}

{\sc Proof} of Lemma \ref{degABN}. Let $F^n H^*$ be the filtration of the limit of the spectral sequence of the first quarter. It consists of a sequence of inclusions
\begin{multline*}
H^n = F^0 H^n \underset{E^{0n}_{\infty}}{\supseteq} 
F^1H^n  \supseteq \cdots \supseteq 
F^{p}H^n \underset{E^{p\,n-p}_{\infty}}{\supseteq} F^{p+1}H^n 
\supseteq \cdots
\\
\cdots \supseteq F^{n-1}H^n \underset{E^{n-1\,1}_{\infty}}{\supseteq} 
F^n H^n  \underset{E^{n0}_{\infty}}{\supseteq} F^{n+1}H^n=0.
\end{multline*}
In this diagram, each symbol under the inclusion denotes their quotient-object 
$E^{p\,n-p}_{\infty}= F^p H^n/ F^{p+1} H^n$ for $0\leq p\leq n$.
 
  For each $n\geq 0$, there exists a decomposition of the edge morphism $\gamma: E^{n0}_2\to H^n$ into the composition
  $$
  E^{n0}_2 \twoheadrightarrow E^{n0}_{\infty}= F^n H^n \subseteq H^n.
  $$

  A morphism $\gamma$ is an isomorphism for all $n$ in the range $0\leq n\leq N$ if and only if $E^{n0}_2\to E^{n0}_{\infty}$ are isomorphisms for $0\leq n\leq N$ and $F^nH^n= H^n$.

A morphism $\gamma$ is a monomorphism for $n=N+1$ if and only if $E^{n0}_2 \twoheadrightarrow E^{n0}_{\infty}$ is a monomorphism, and hence an isomorphism.

Therefore, the $(A_N)$ condition is equivalent to the following two conditions:
\begin{itemize}
\item $E^{n 0}_2\to E^{n 0}_{\infty}$ is an isomorphism as $0\leq n\leq N+1$;
\item $E^{0n}_{\infty}=0, \ldots , E^{n-1\,1}_{\infty}=0$, for $1\leq n\leq N$.
\end{itemize}

Let us prove the implication $(B_N)\Rightarrow (A_N)$ under the conditions of the lemma.
The execution of $(B_N)$ means that $E^{0,n}_2=0$ is true for all $n$ such that $1\leq n\leq N$.
By the condition of the lemma, in this case $E^{pq}_2=0$ is true for all $p\geq 0$ and $1\leq q\leq N$.
Hence $H^n=F^0 H^n= \ldots F^n H^n= E^{n,0}_{\infty}= E^{n, 0}_2$, and for $n= N+1$ we get $E^{N+1,0}_{\infty}\subseteq H^{N+1}$.
(May be proved by dual to \cite[Exercise XI.3.2]{mac1963}.)

Let us prove the implication $(A_N)\Rightarrow (B_N)$ by induction on $N\geq 1$.
Let's check for $N=1$. Let $E^{1,0}\to H^1$ be an isomorphism and $E^{2,0}\to H^2$ be a monomorphism. Consider the exact sequence of lower powers
$$
0 \to E^{1,0}_2\to H^1\to E^{0,1}_2 \to E^{2,0}_2 \to H^2
$$
(see the dual exact sequence in \cite[Exercise XI.3.2]{mac1963}).
The conditions $(A_N)$ for $N=1$, together with this exact sequence, lead to the isomorphism $E^{0,1}_2=0$, which means the fulfillment of $(B_1)$.

Let the implication be true for $N$, i.e. $(A_N)$ implies $E^{0,n}_2= 0$ for $n= 1, \ldots, N$. It suffices to prove that $(A_{N+1})$ implies that $E^{0,n}_2=0$ for $n=1, \ldots, N+1$.
Consider the exact sequence
$$
  E^{0\,N+1}_{\infty} \to E^{0\,N+1}_2 \to E^{N+2\,0}_2\to E^{N+2\ ,0}_{\infty} \to 0
$$
Property $(A_{N+1})$ implies that $E^{0\,N+1}_{\infty}= 0$ and $E^{N+2\,0}_2\to E ^{N+2\,0}_{\infty}$ is an isomorphism, whence $ E^{0\,N+1}_2 =0$.
\hfill$\Box$

\begin{proposition}\label{condVerdier}
For any $N\geq 1$, the following properties of the functor $f: \mC\to \mD$ are equivalent:
\begin{enumerate}
\item For each diagram $\mG: \mD\to \Ab$ the homomorphisms 
$H^n(\mD, \mG)\to H^n(\mC, \mG\circ f)$ are isomorphisms for all $0\leq n\leq N$ and monomorphisms for $n=N+1$.
\item For each $d\in \Ob\mD$ the homology group $H_n(B(f\downarrow d))=0$ for all $1\leq n\leq N$ and the category $f\downarrow d$ is connected and non-empty .
\end{enumerate}
\end{proposition}
\begin{corollary}(Verdier criterion)
Canonical homomorphisms
$$
H^n(\mD, \mG)\to H^n(\mC, \mG\circ f)
$$
are isomorphisms for all $\mG\in \Ab^{\mD}$ and $n\geq 0$ if and only if for each $d\in \Ob\mD$ of the group $H_n(f\downarrow d)$ are equal to $0$ for all $n>0$ and are isomorphic to $\ZZ$ for $n=0$.
\end{corollary}

In particular, if the functor $f$ has a right adjoint, then for every $d\in \Ob\mD$ the category $f\downarrow d$ has a final object, and the Verdier criterion gives isomorphisms $H^n(\mD, \mG ) \xrightarrow{\cong} H^n(\mC, \mG\circ f)$ for all $\mG\in \Ab^{\mD}$ and $n\geq 0$.

\section{Baues-Wirsching cohomology}

In this section we recall the definition of the category of factorizations for the small category and the definition of the Baues-Wirsching cohomology groups $H^n_{BW}(\mC, \mG)$ for the small category $\mC$ with coefficients in the natural system $\mG: \fF\mC\to \Ab$, and consider an example of calculating these cohomology groups.
After that, we give an example of a left adjoint functor $f: \mC\to \mD$ between small categories and a natural system $\mG$ on $\mD$ such that $H^1_{BW}(\mD, \mG) \to H^1_{BW}(\mD, f^*\mG)$ is not an isomorphism.
Then we give a characterization of functors $f: \mC\to \mD$ such that for every natural system $\mG$ on $\mD$ and every $n\geq 0$ the canonical homomorphism $H^n_{BW} (\mD, \mG)\to H^n_{BW}(\mD, f^*\mG)$ is an isomorphism.

\subsection{Category of factorizations}

Let $\mC$ be a small category. The factorization category $\fF\mC$ \cite{bau1985} is defined as follows.

Its objects are the morphisms $f\in \Mor\mC$. Its morphisms $f\to g$ are given by the pairs $(\alpha,\beta)\in \Mor\mC^{op}\times \Mor\mC$ making the diagram commutative in the category $\mC$:
\begin{equation}\label{sqmor}
\xymatrix{
\circ \ar[r]^{\beta}& \circ \\
\circ \ar[u]^f & \circ \ar[l]^{\alpha} \ar[u]_g
}
\end{equation}
We will denote morphisms by $f \xrightarrow{(\alpha, \beta)}g$.
The identity morphism of the object $f\in \Ob\fF\mC$ is equal to $f \xrightarrow{(id_{\dom f}, id_{\cod f})} f$.
The morphism composition $f \xrightarrow{(\alpha, \beta)}g \xrightarrow{(\alpha', \beta')}h$ is defined as $f \xrightarrow{(\alpha\alpha', \beta' \beta )}h$.

A natural system on a small category $\mC$ is a diagram of Abelian groups $\mG: \fF\mC\to \Ab$.
For this natural system, as for a diagram on the category of factorizations, the cohomology groups $H^n(\fF\mC, \mG)$ are defined for all $n\geq 0$.

Let us compare the nerves of the categories $\mC$ and $\fF\mC$.
Consider the functor $\cod: \fF\mC\to \mC$, which assigns to each $\alpha\in \Ob\fF\mC$ its end, and to the morphism $\alpha\xrightarrow{(u,v)}\beta$ morphism $v: \cod\alpha \to \cod\beta$.

\begin{proposition}\label{codqui}
Let $\mC$ be a small category.
For each $c\in \Ob\mC$ the category $\cod\downarrow c$ contains a reflective subcategory, in the sense of \cite[Definition 3.5.2]{bor1994}, isomorphic to the category $(\mC\downarrow c)^{op }$. Hence, the nerve of the category $\cod\downarrow c$ is contractible. In particular, the nerve map $B(\cod): B(\fF\mC)\to B\mC$ is a weak equivalence.
\end{proposition}
Weak equivalence follows from Quillen's theorem \cite[Theorem A]{qui1973}.

 \subsection{Definition and examples of Baues-Wirsching cohomology}

Let $\mC$ be a small category and $\mG: \fF\mC\to \Ab$ be a natural system on $\mC$.

Consider the complex $(F^*(\mC,\mG), d^*)$ of Abelian groups
$$
F^n(\mC, \mG)= \prod_{c_0\xrightarrow{\alpha_1} \cdots \xrightarrow{\alpha_n} c_n} \mG(\alpha_n\cdots\alpha_1), \quad n\geq 0 ,
$$
whose elements are defined as functions $\varphi: B_n\mC\to \bigcup_{\alpha\in \Mor\mC} \mG(\alpha)$ taking the values $\varphi(c_0\xrightarrow{\alpha_1} c_1 \to \cdots \to c_{n-1}\xrightarrow{\alpha_n}c_n)\in \mG(\alpha_n\cdots \alpha_1)$.
The differentials are $d^n= \sum^n_{i=0}(-1)^i d^n_i$, where $d^n_i: F^n(\mC, \mG)\to F^{n+1 }(\mC, \mG)$ are defined by the formulas \\
$d^n_i\varphi(c_0\xrightarrow{\alpha_1} c_1\to \ldots \to c_{n}\xrightarrow{\alpha_{n+1}}c_{n+1})=$
 $$
 \begin{cases}
 \mG(\alpha_1, 1_{c_0}) (\varphi(c_1\xrightarrow{\alpha_2} c_2\to \ldots \to c_{n}
  \xrightarrow{\alpha_{n+1}} c_{n+1})), & i = 0,\\
\varphi(c_0\xrightarrow{\alpha_1} c_1\to \ldots \to \widehat{c_i}\to \ldots \to c_{n}
 \xrightarrow{\alpha_{n+1}} c_{n+1}), & 1\leq i \leq n,\\
 \mG(1_{c_{n+1}}, \alpha_{n+1})
  (\varphi(c_0\xrightarrow{\alpha_1} c_1\to  \ldots \xrightarrow{\alpha_{n}} c_{n})), 
  & i = n+1. 
 \end{cases}
 $$

\
\begin{definition}
The cohomology of the complex $(F^*(\mC, \mG), d^*)$ is called the Baues-Wirsching cohomology $H^n(\mC, \mG)$ of the small category $\mC$ with coefficients in the natural system $\mG$, $n\geq 0$.
\end{definition}

Baues and Wirsching proved that the cohomology of $H^n_{BW}(\mC, \mG)$ is isomorphic to $H^n(\fF\mC, \mG)$ for all $n\geq 0$ \cite[Theorem 4.4 ]{bau1985}.

In Example \ref{mon1}, we have established that the cohomology groups of an idempotent monoid are equal to 0 in all positive dimensions.

Consider the Baues-Wirsching cohomology of this monoid.

\begin{example}
The idempotent monoid $\mE$ consists of two elements: $1$ and $e$, with $e^2=e$. The factorization category $\fF\mE$ consists of two objects $1$ and $e$ and eight morphisms:

\centerline{$(1, 1): 1\to 1, (1,e): 1\to e, (e,1): 1\to e, (e,e): 1\to e,$}
\centerline{$(1, 1): e\to e, (1,e): e\to e, (e,1): e\to e, (e,e): e\to e.$}
Details of the calculation can be found in \cite[Example 2.4]{X1997}.
For the natural system $\cA: \fF\mE\to \Ab$ the cohomology $H_{BW}^n(\mE, \cA)$ will be isomorphic to the cohomology groups of the complex
$$
0\to \cA(1)\xrightarrow{d^0} \cA(e)\xrightarrow{d^1} \cA(e)\xrightarrow{d^2} \cA(e)\xrightarrow{d^3 }\to\ldots,
$$
with differentials
$$
d^n(\varphi)=
\begin{cases}
F(e,1)\varphi - F(1,e)\varphi \text{ for even $n$,}\\
F(e,1)\varphi - \varphi +F(1,e)\varphi \text{ for odd $n$.}
\end{cases}
$$
This complex is exact in dimensions $n \geq 2$.
If $\cA(1)= 0$ and $\cA(e)= \ZZ\oplus \ZZ$, then $H^1_{BW}(\mE, \cA)= \ZZ\oplus\ZZ$, and $H^n_{BW}(\mE, \cA)= 0$ for all $n\not=0$.
\end{example}

 The following example shows that the left adjoint functor $f: \mC\to \mD$ may not 
 preserve cohomology with coefficients in the natural system.
 
  \begin{example}\label{BWNotIso}
  Let $\mC= [0]$ be the category consisting of the unique object $0$ and the unique morphism $id_0$ equal to the identity. Consider the category $\mD=[1]$, which consists of two objects $0$ and $1$, a morphism $\gamma: 0\to 1$, and two identical morphisms $id_0$ and $id_1$.
  The factorization category $\fF\mC$ will consist of a single object and a single morphism (equal to the identity). The category $\fF\mD$ will have three objects and two morphisms besides the identity, it is isomorphic to the category $V$ from the example \ref{pos1} and looks like this:
  $$
  id_0 \to \gamma \leftarrow id_1.
  $$
  Consider a natural system $D: \fF\mD\to \Ab$ that takes the values $D(\gamma)=\ZZ$ and the remaining $D(id_0)= D(id_1)=0$ are equal to zero.
  As in the example \ref{pos1}, the first cohomology group $H^1_{BW}(\mD, D)= H^1(\fF\mD, D)$ is isomorphic to $\ZZ$. The rest, for $n\not=1$, are equal to $H^n_{BW}(\mD, D)=0$. For any natural system $E$ on $\mC=[0]$, the cohomology groups are equal to $H^0_{BW}(\mC, E)= H^0(\fF\mC, E)=E(0) $ and $H^n_{BW}(\fF\mC, E)=0$ for all $n>0$.
   Hence $H^1_{BW}(\mD, D)\to H^1_{BW}(\mC, f^*D)$ is not an isomorphism.

Consider the functor $f: \mC\to \mD$ taking the values $f(0)=0$. It will have a right adjoint
functor. Therefore, a left adjoint functor may not preserve the Baues-Wirsching cohomology groups.
Note that we can take $f(0)= 1$. Then $f$ will have a left adjoint. And in this case $H^1_{BW}(\mD, D)=\ZZ$, and $H^1_{BW}(\mC, f^*D)=0$ are not isomorphic. The right adjoint functor, like the left one, does not preserve the Baues-Wirsching cohomology.
  \end{example}
 
 \subsection{Preservation of the Baues-Wirsching cohomology}
 
 Here we study the left fibres of the functor between small categories and describe a criterion for the preservation of the Bowes-Wirsching cohomology when passing to the inverse image.
 
  Let $f: \mC\to \mD$ be a functor between small categories.
For $\alpha\in \Mor\mD$, denote by $f\langle\alpha\rangle$ the category whose objects are pairs of morphisms $d \xrightarrow{\beta} f(c) \xrightarrow{\beta'} d' $ such that $\beta'\beta=\alpha$. Morphisms between $d \xrightarrow{\beta} f(c) \xrightarrow{\beta'} d'$ and $d \xrightarrow{\gamma} f(c') \xrightarrow{\gamma'} d'$ are given with morphisms $\nu: c\to c'$ making the following diagram commutative
 $$
 \xymatrix{
 & f(c) \ar[dd]^{f(\nu)} \ar[rd]^{\beta'}\\
 d \ar[ru]^{\beta} \ar[rd]_{\gamma} && d' \\
 & f(c') \ar[ru]_{\gamma'}
 }
 $$

 \begin{lemma}\label{isocats}
  For every functor $f: \mC\to \mD$ and for every $\alpha\in \Mor\mD$, the category $\fF(f)\downarrow\alpha$ is isomorphic to $\fF(f\langle\alpha\rangle)$.
  \end{lemma}
  {\sc Proof.} We construct a functor that assigns to each object $(f(\beta) \xrightarrow{(u,v)}\alpha) \in \Ob(\fF f \downarrow\alpha)$ an object $(u, vf\beta)\xrightarrow{\beta} (f(\beta)u, v)$ of the category $\fF(f\langle\alpha\rangle)$.
The action of a functor on morphisms $(u,v)\xrightarrow{(\gamma, \gamma')} (u',v')$ of the category $\fF(f)\downarrow \alpha$ can be shown using the following commutative diagrams
$$
\xymatrix{
f(\beta)\ar[d]_{(\gamma,\gamma')} \ar[r]^{(u,v)} & \alpha\\
f(\beta') \ar[ru]_{(u',v')}
}
\qquad
\xymatrix{
\ar@{}[d]_{\mapsto} \\
&
}
\xymatrix{
(u, vf(\beta)) \ar[r]^{\beta} & (f(\beta)u, v) \ar[d]^{\gamma'} \\
(u', v'f(\beta')) \ar[u]_{\gamma} \ar[r]_{\beta'} & (f(\beta')u', v')
}
$$
The diagram on the left shows a morphism in the category $\fF f \downarrow \alpha$. On the right is the corresponding morphism from $\fF(f\langle\alpha\rangle)$.
If we add the morphism $(u',v')\xrightarrow{(\zeta, \zeta')} (u'',v'')$ on the left, then on the right we get one more rectangle located from below.
On the right we get two commutative rectangles, proving the commutativity of the enclosing rectangle.
Thus, the composition of mor\-phisms of the category $\fF f\downarrow \alpha$ goes over to composition, and hence the constructed mapping is functorial. It is clear that this mapping is bijective.
Therefore, it will be an isomorphism.
\hfill$\Box$
 
  \begin{theorem}
  Let $f: \mC\to \mD$ be a functor between small categories.
  For every $N\geq 1$ the following properties of this functor are equivalent:
  \begin{itemize}
  \item For any diagram $\mG: \fF\mD \to \Ab$ the canonical homomorphisms $H^n_{BW}(\mD, \mG) \to H^n_{BW}(\mC, \mG\circ f)$ are isomorphisms for all $0\leq n\leq N$ and monomorphisms for $n=N+1$.
  \item The homology groups $H_n(f\langle\alpha\rangle)$ are zero for all $\alpha\in \Mor\mD$ for $1\leq n\leq N$ and the category $f\langle\alpha\rangle$ is connected and non-empty.
  \end{itemize}
  \end{theorem}
  {\sc Proof.} According to \cite[Theorem 4.4]{bau1985}, the cohomology $H^n_{BW}(\mD, \mG)$ is naturally isomorphic to the cohomology $H^n(\fF\mD, \mG)$ .
  We use Proposition \ref{condVerdier} for the functor $\fF{f}: \fF\mC\to \fF\mD$ and the diagram $\mG: \fF\mD\to \Ab$. We get the equivalence of statements
\begin{enumerate}
\item $H^n(\fF\mD, \mG)\to H^n(\fF\mC, \mG\circ \fF{f})$ are isomorphisms for all $0\leq n\leq N$ and monomorphisms for $n=N+1$.
\item For each $\alpha\in \Ob\fF\mD$, the homology groups $H_n(B(\fF{f}\downarrow \alpha))=0$ for all $1\leq n\leq N$ and the category $ \fF{f}\downarrow \alpha$ is connected and non-empty.
\end{enumerate}
The category $\fF{f}\downarrow \alpha$ is isomorphic to $\fF{f\langle\alpha\rangle}$ by Lemma \ref{isocats}. According to Proposition \ref{codqui}, for any category its nerve is weakly equivalent to its factorization category, which implies the isomorphism $H_n(B(\fF{f}\downarrow\alpha))\cong H_n(B(f\langle\alpha\rangle))$ for all $n\geq 0$.

Thus, the proved theorem is obtained from the Proposition \ref{condVerdier} by substituting the functor $\fF{f}$ instead of the functor $f$.
  \hfill$\Box$

 \section{Hochschild-Mitchell cohomology}
 
  In this section, we recall the definition of Hochschild-Mitchell cohomology and the proof of the theorem that these cohomology can be computed as Baues-Wirsching cohomology.
  Then we give an example of the left adjoint functor $f: \mC\to \mD$ for which the groups $H^1_{HM}(\mD, \mG)$ and $H^1_{HM}(\mC, \mG\circ(f^{op}\times f))$ are not isomorphic.
  We also prove a theorem on the characterization of Hochschild-Mitchell cohomology isomorphism-inducing functors.

  \subsection{Definition and examples of Hochschild-Mitchell cohomology}
 
  For an arbitrary small category $\mC$ we denote by $\ZZ\mC: \mC^{op}\times \mC\to \Ab$ the composition of the morphism bifunctor $\mC^{op}\times \mC\xrightarrow{\mC (-,=)} \Set$
  and the functor $\Set \to \Ab$, which assigns to each set a free Abelian group generated by this set and a map - a homomorphism of free Abelian groups that extends this map.

  \begin{definition}\cite[Definition 8.1]{bau1985}
  Let $\mC$ be a small category and $\mG: \mC^{op}\times \mC\to \Ab$ be a diagram of abelian groups. The Hochschild-Mitchell cohomology of the category $\mC$ with coefficients in $\mG$ is defined by the formula
  $$
 H^n_{HM}(\mC, \mG):= Ext^n_{\Ab^{\mC^{op}\times \mC}}(\ZZ\mC, \mG), 
 n\geq 0.
  $$
  \end{definition}
Mitchell gave this definition for the pre-additive category \cite[\S12]{mit1972}.

An example of computing the Hochschild-Mitchel cohomology for an idempotent monoid is given in \cite[Example 2.4]{X1997}.
 
  Consider the functor $(\dom, \cod): \fF\mC \to \mC^{op}\times \mC$, which assigns to each $\alpha\in \Ob\fF\mC$ an object $(\dom\alpha , \cod\alpha)\in \Ob(\mC^{op}\times\mC)$, and a pair of morphisms $(f,g): \alpha\to \beta$ - a pair of the same morphisms $(f, g)$ in the category $\mC^{op}\times \mC$. The objects of the left fibre $(\cod, \dom)\downarrow (a,b)$ we will define as triples of morphisms $(\alpha, (x,y))$, as shown in the following figure
 $$
 \xymatrix{
  \cod\alpha \ar[rr]^y && b\\ 
  \dom\alpha\ar[u]^{\alpha} && \ar[ll]^{x} a
 }
$$
 
  \begin{lemma}\cite[Lemma 2.4]{X1997}
  For any object $(a,b)\in \mC^{op}\times \mC$ the connected component of the left fibre
   $(\dom, \cod)\downarrow (a,b)$ containing the object $(\alpha, (x, y))$ 
   has a final object $(y\circ \alpha\circ x, (1_a, 1_b))$.
  \end{lemma}
 
  This lemma implies that $Lan^{(\dom, \cod)}\Delta_{\fF\mC}\ZZ=\ZZ\mC$ and $Lan^{(\dom, \cod)}_q\Delta_ {\fF\mC}\ZZ=0$ for all $q>0$.

Let's apply Proposition \ref{spobe} to the functor $(\dom, \cod): \fF\mC\to \mC^{op}\times \mC$.
The resulting spectral sequence degenerates and leads to isomorphisms:

  \begin{corollary}\cite[Proposition 8.5(A)]{bau1985}
  For any diagram $\mG: \mC^{op}\times \mC\to \Ab$ and $n\geq 0$ the formula $H^n_{HM}(\mC, \mG) \cong H^n_ {BW}(\mC, \mG\circ(\dom, \cod))$.
  \end{corollary}
 
  As the following example shows, the left adjoint functor does not preserve the Hochschild-Mitchel cohomology.
 
 \begin{example}
  Consider the category $\mD=[1]$, as in Example \ref{BWNotIso}.
  The category $\mD^{op}\times \mD$ consists of the following objects and morphisms (except identical ones)
 $$
 \xymatrix{
 (1,0)\ar[d]_{(\gamma^{op}, id_0)} \ar[r]^{(id_1,\gamma)} 
 \ar[rd]|-{(\gamma^{op}, \gamma)}& 
 (1,1) \ar[d]^{(\gamma^{op}, id_1)}
 \\
 (0,0) \ar[r]_{(id_0,\gamma)} & (0,1)
 }
 $$
  The functor $(\dom, \cod)$ maps the category $\fF\mD$ to a full subcategory
  of the category $\mD^{op}\times\mD$ that has objects
 $$
 \xymatrix{
  & (1, 1)\ar[d]\\
  (0,0) \ar[r] & (0,1)
 }
 $$
 Consider the diagram of Abelian groups $\mG: \mD^{op}\times \mD\to \Ab$, which takes values on objects equal to the zero Abelian group except for $\mG(0,1)=\ZZ$.
  The value on the identity morphism of the object $(0,1)$ is equal to the identity homomorphism. On all other morphisms, the values of the bimodule $\mD$ are equal to the zero homomorphisms.
The groups $H^n_{HM}(\mD, \mG)$ are isomorphic to the groups
$H^n_{BW}(\mD, \mG\circ(\dom,\cod))$ computed in Example \ref{BWNotIso} 
where the diagram $\mG\circ(\dom, \cod)$ denoted by $D$.
Means
$$
H^0_{HM}(\mD, \mG)= 0, ~H^1_{HM}(\mD,\mG)= \ZZ \text{ and } H^n_{HM}(\mD, \mG)=0
\text{ for all } n\geq 2.
$$

 Let $\mC=[0]$, then the categories $\mC^{op}\times \mC$ and $\fF\mC$ are isomorphic to $[0]$.
If the functor $f: \mC\to \mD$ is defined as $f(0)=0$, then it is left adjoint, and in the case $f(0)=1$ it is right adjoint.
For any functor $f: \mC\to \mD$ of the group $H^n_{HM}(\mC, \mG\circ(f^{op}\times f))= H^n_{BW}(\mC , \mG\circ(f^{op}\times f)\circ(\dom, \cod))$ are equal to $0$ for all $n>0$. Hence $H^1_{HM}(\mD, \mG)\to H^1_{HM}(\mC, \mG\circ(f^{op}\times f^{op}))$ is not an isomorphism .
 \end{example} 
 
 \subsection{Preservation of the Hochschild-Mitchell cohomology}
 
  Using the Oberst spectral sequence (Proposition \ref{oberst1967}), we study the invariance conditions for the Hochschild-Mitchel cohomology of bimodules under the transition to the inverse image
 $$
 E^{pq}_2 = Ext^p(Lan^{f^{op}\times f}_q\ZZ\mC, \mG) \Rightarrow 
 Ext^n(\ZZ\mC, \mG\circ(f^{op}\times f))
 $$
 
    First we construct canonical homomorphisms
  \begin{equation}\label{canon3}
  H^0_{HM}(\mD, \mG)\to H^0_{HM}(\mC, \mG\circ(f^{op}\times f))
  \end{equation}
   There is a natural transformation
   $\ZZ\mC\to \ZZ\mD\circ(f^{op}\times f)$. The left adjoint functor to $(-)\circ{f^{op}\times f}$ induces a morphism $Lan^{f^{op}\times f}\ZZ\mC\to \ZZ\mD$. It defines for each $\mG\in \Ab^{\mD^{op}\times \mD}$ a canonical homomorphism  
  $$
  Hom(\ZZ\mD, \mG)\to Hom(Lan^{f^{op}\times f}\ZZ\mC, \mG)= Hom(\ZZ\mC,
\mG\circ (f^{op}\times f) ) 
  $$
  that defines a canonical morphism (\ref{canon3}).
  
  \begin{lemma}\label{HM0}
   A homomorphism (\ref{canon3}) is an isomorphism for all $\mG\in \Ab^{\mD^{op}\times\mD}$ if and only if the natural transformation $Lan^{f^{op}\times f}\ZZ\mC\to
   \ZZ\mD$ is an isomorphism.
  \end{lemma}
  
 \begin{theorem}
 Let $f: \mC\to \mD$ be a functor between small categories and $N\geq 1$ be a positive integer.
  Then the following two properties of this functor are equivalent
  \begin{enumerate}
  \item $Lan^{f^{op}\times f}\ZZ\mC \to \ZZ\mD$ is an isomorphism and for all $1\leq n\leq N$ all values of the diagram $Lan^{f^{op }\times f}_n\ZZ\mC$ are equal to the zero group.
  \item For all $0\leq n\leq N$ the homomorphisms $H^n_{HM}(\mD, \mG)\to H^n_{HM}(\mC, \mG\circ(f^{op}\times f))$ are isomorphisms, and for $n=N+1$ they are monomorphisms, for all bimodules $\mG: \mD^{op}\times \mD\to \Ab$.
  \end{enumerate}
 \end{theorem}
The proof follows from Lemma \ref{HM0} and Proposition \ref{oberst1967}.

 \section{Thomason cohomology}
 
  We recall the definition of the Thomason cohomology $H^n_{T}(\mC, \mG)$ and give an example showing that a left adjoint functor may not preserve the Thomason cohomology.
  Then we prove a theorem on the characterization of functors between small categories that preserve the Thomason cohomology under the inverse image.
 
  \subsection{Thomason cohomology definition}
 
  Let $\mC$ be a small category.
  Consider the category $\Delta\downarrow \mC$ whose set of objects consists of functors $\sigma: [m]\to \mC $, $m\geq 0$, and morphisms between objects $\sigma\xrightarrow{\alpha} \sigma'$ are commutative triangles
 $$
 \xymatrix{
 [m] \ar[rd]_{\alpha} \ar[rr]^{\sigma} && \mC \\
 & [n] \ar[ru]_{\sigma'} 
 }
 $$

 \begin{definition}\cite{gal2013}
  The Thomason cohomology $H^n_T(\mC, \mG)$ of the small category $\mC$ with coefficients in the diagram $\mG: \Delta\downarrow \mC\to \Ab$ are the cohomology groups of the cochain complex consisting of the products 
  $C^n= \prod\limits_{\sigma\in \mC^{[n]}}G(\sigma)$ and differentials
  $d^n=\sum^{n+1}_{i=0}(-1)^i d^n_i: C^n\to C^{n+1}$ given as
  $$
  (d^n_i\psi)(\sigma')=
  \mG(\sigma\xrightarrow{\partial^i_{n+1}}\sigma')
  (\psi(\sigma)).
  $$
   \end{definition}

   It is known that the Thomason cohomology with coefficients in $\mG: \Delta\downarrow \mC\to \Ab$ is isomorphic
    $H^n(\Delta\downarrow \mC, \mG):= \Lim^n_{\Delta\downarrow \mC,}\mG$.
    It follows from the dual statement to
\cite[Appendix 2, Proposition 4.2]{gab1967}, for the case $X=B\mC$.
 
Let us give a counterexample showing that the left adjoint functor does not preserve the Thomason cohomology under the inverse image.
 
  Note first that for any functor $G: \fF\mC \to \Ab$ the cohomology isomorphism $H^n(\mC, G)\xrightarrow{\cong} \Lim^n_{\Delta/\mC} (G\circ\delta_{\mC})$
  for all $n\geq 0$. Here $\delta: \Delta/\mC\to \fF\mC$ is a functor assigning each sequence of morphisms of the diagram $\sigma: [n]\to \mC$ to their composition.
  This is proved in \cite[Theorem 4.4]{bau1985}.
  It can be proved with the help of the following assertion, which we will need more.

 \begin{lemma}\label{delqui}
  For any category $\mC$ and $\alpha\in \Ob\fF\mC$, the nerve of the category $\delta \downarrow \alpha $ is weakly equivalent to a point.
  \end{lemma}
  {\sc Proof.} Let $c=\dom\alpha$ be the beginning of the morphism $\alpha$. The category $\delta\downarrow\alpha$ is isomorphic to $\Delta\downarrow((c\downarrow\mC)\downarrow \alpha)$.
For an arbitrary small category $\mD$, the functor $\partial: \Delta\downarrow \mD \to \mD$, $\partial(c_0\to c_1\to \ldots \to c_n)=c_n$, has contractible left fibres, whence the nerve of the category $\Delta\downarrow \mD$ is weakly equivalent to $B\mD$. Hence the nerve
 $\delta\downarrow \alpha$ is weakly equivalent to the nerve of the category
  $(c\downarrow\mC)\downarrow \alpha$ having an initial (and final) object. Hence the nerve $\delta\downarrow \alpha$ is weakly equivalent to a point.
\hfill$\Box$
 
 \begin{example}
 Consider the natural system $D: \fF\mD\to \Ab$ and the following commutative diagram
  $$
  \xymatrix{
  \ar[d]_{\delta_{\mC}}\Delta/\mC \ar[rr]^{\Delta/f} &&
  \Delta/\mD\ar[d]^{\delta_{\mD}}\\
  \ar[rd]_{D\circ\fF f} \fF\mC \ar[rr]^{\fF f} && \fF\mD \ar[ld]^{D} \\
  & \Ab
  }
  $$
Take the categories $\mC$, $\mD$, the functor $f$, and the natural system $D$ from 
Example \ref{BWNotIso}.
Then the groups $H^1_{BW}(\mD, D)$ and $H^1_{BW}(\mC, D\circ\fF f)$ are not isomorphic. For all $n\geq 0$ there are isomorphisms $\Lim^n_{\fF\mD} D \cong \Lim^n_{\Delta/\mD} D\delta_{\mD}$ and
$$
\Lim^n_{\fF\mC} D\fF f
\cong \Lim^n_{\Delta/\mC} D(\fF f)\delta_{\mC}=
  \Lim^n_{\Delta/\mD} D\delta_{\mD}(\Delta/f)
$$
where the last isomorphism is an equality since the outer pentagon is commutative. Hence the Thomason cohomology groups $H^1_T(\mD, D\delta_{\mD})\cong \Lim^1_{\Delta/\mD} D\delta_{\mD}$ and $H^1_T(\mD , D\delta_{\mD}(\Delta/f))\cong \Lim^1_{\Delta/\mD} D\delta_{\mD}(\Delta/f)$ are not isomorphic.
 \end{example}
 
 \subsection{Preservation of the Thomason cohomology}
  
Let us describe the necessary and sufficient conditions for the functor $f: \mC\to \mD$ under which for each diagram of Abelian groups $\mG: \Delta\downarrow \mD\to \Ab$ the homomorphisms $H^n_T(\mD, \mG)\to H^n_T(\mC, \mG\circ (\Delta\downarrow f))$ are invertible for all $n\geq 0$.
 
  Let $f: \mC\to \mD$ be a functor between small categories.
 For an arbitrary functor $\sigma: [n]\to \mC$, let $\overleftarrow{f}({\sigma})$ denote the category that makes the next square Cartesian (pullback) in $\Cat$ with sides $f$ 
 and $\sigma$:
  $$
  \xymatrix{
  \mC \ar[r]^f & \mD\\
  \overleftarrow{f}({\sigma})\ar[u] \ar[r] & [n]\ar[u]_{\sigma}
  }
  $$
Due to the universal property of this square, the formula $(\Delta\downarrow f)\downarrow \sigma\cong \Delta\downarrow \overleftarrow{f}(\sigma)$ will be true.
The nerves of the categories $\Delta\downarrow \overleftarrow{f}(\sigma)$ and $\overleftarrow{f}(\sigma)$ are weakly equivalent.
(This follows from the well-known assertion that for an arbitrary small category $\mC$ there exists a weak equivalence $\Delta\downarrow \mC\to \mC$, which can be obtained from Lemma \ref{delqui} and Proposition \ref{codqui}.)
Therefore, the following is true.

\begin{lemma}\label{comqui}
The nerve of the left fibre for the functor $\Delta\downarrow f: \Delta\downarrow \mC \to \Delta\downarrow \mD$ for each $\sigma\in \Ob(\Delta\downarrow \mD)$ is weakly equivalent to the nerve of the category $\overleftarrow{f}(\sigma)$.
\end{lemma}
 
 Since $H^n_{T}(\mC, \mG)\cong \Lim^n_{\Delta\downarrow\mC}\mG$, we can use the following criterion to characterize Thomason cohomology-preserving functors.
 
\begin{theorem}
For any functor between the small categories $f: \mC\to \mD$ and $N\geq 1$ the following two properties are equivalent
\begin{itemize}
\item Canonical homomorphisms
$$
H^n_T(\mD, \mG)\to H^n_T(\mC, \mG\circ (\Delta\downarrow f))
$$
are isomorphisms for all $\mG\in \Ab^{\Delta\downarrow\mD}$ and $0\leq n\leq N$, and monomorphisms for $n=N+1$.
\item For each $\sigma\in \Ob(\Delta\downarrow\mD)$ the groups $H_n(\overleftarrow{f}({\sigma}))$ are equal to $0$ for all $1\leq n\leq N$ and are isomorphic to $\ZZ$ for $n=0$.
\end{itemize}
\end{theorem}

This follows from 
Lemma \ref{degABN}, Lemma \ref{comqui} and Proposition
\ref{spobe}.

\section{Compatible cohomology of small categories}

We generalize the considered cohomology of small categories. For generalized cohomology, there will also be criteria for preserving the inverse image cohomology.

\subsection{Definition of compatible cohomologies}

We assume that for each small category $\mC$ a small category $\FF\mC$ and a diagram of Abelian groups $\DD\mC: \FF\mC\to \Ab$ are defined, and this pair functorially depends on $\mC $. The term functoriality will be explained later.

\begin{definition}
The compatible $n$-cohomology of the small category $\mC$ with coefficients
 in the natural system $\mG: \FF\mC\to \Ab$ is the Abelian group 
 $H^n_{\DF}(\mC, \mG)= Ext ^n_{\Ab^{\FF\mC}}(\DD\mC, \mG)$.
\end{definition}

\begin{example}
The Hochschild-Mitchel cohomologies are compatible.
For them $\FF\mC=\mC^{op}\times \mC$, $\DD\mC= \ZZ\mC: \mC^{op}\times \mC\to \Ab$.
For any bimodule $\mG: \mC^{op}\times \mC\to \Ab$ the compatible cohomologies are $H^n_{\DF}(\mC, \mG)=
Ext^{n}_{\Ab^{\mC^{op}\times \mC}}(\ZZ\mC, \mG)$.
Hence, they are isomorphic to the Hochschild-Mitchel cohomology.
\end{example}

Consider the Grothendieck construction $\Cat^{op}\int \Ab^{\DD(-)}$ for the functor $\Ab^{\FF(-)}: \Cat^{op}\to \CAT$, matching to each small category $\mC$ the diagram category $\Ab^{\FF\mC}$ and transferring each morphism $f: \mC\to \mD$ of the category $\Cat$ to the functor $\Ab^{\FF( f)}: \Ab^{\mD}\to \Ab^{\mC}$ acting as $\mG\mapsto \mG\circ \FF(f)$.

Let us explain what the functoriality of the pair $(\FF, \DD)$ means.

For the functor $\Ab^{\FF(-)}: \Cat^{op}\to \CAT$ there exists a contravariant Grothendieck construction, which is defined in \cite{nlabGC}.
Denote it by $\Cat\int^{*}\Ab^{\FF(-)}$. The objects of this category are the pairs $(\mC, \mF)$, consisting of $\mC\in \Cat$ and the diagram $\mF: \FF\mC\to \Ab$, and the morphisms are given as the pairs $(f, \xi)$ consisting of the functor $f: \mC\to \mD$ and the natural transformation $\xi: \mF\to \mG\circ \FF(f)$.
A pair $(\FF, \DD)$ is said to be functorial if the mapping that assigns to each small category $\mC$ an object $\DD\mC: \FF\mC\to \Ab$, and to a morphism $f: \mC\to \mD$ - the pair $(f, \xi: \DD\mC\to \DD\mD\circ \FF(f))$, as the diagram below shows, is a functor $\Cat\to \Cat\int^{*}\Ab^{\FF(-)}$.
$$
\xymatrix{
\FF\mC\ar[rd]_{\DD\mC}^{\xi\nearrow} \ar[rr]^{\FF(f)}  && \ar[ld]^{\DD\mD} \FF\mD\\
& \Ab
}
$$

\begin{proposition}
If a pair of mappings $(\FF, \DD): \Cat \to \Cat\int^*\Ab^{\FF(-)}$ is functorial, then the compatible $n$-cohomology groups $H^ n_{\DF}(\mC, \mG)$, for each $n\geq 0$, define the functor $H^n_{\DF}: \Cat^{op}\int \Ab^{\FF(- )}\to \Ab$.
\end{proposition}

\subsection{Preservation of compatible cohomology}

We consider the functors $(\FF, \DD): \Cat\to \Cat\int^* \Ab^{\FF(-)}$ described above.

\begin{proposition}\label{spscog}
For any functor $(\FF, \DD): \Cat\to \Cat\int^* \Ab^{\FF(-)}$ and functor $f: \mC\to \mD$ there exists a spectral sequence of the first quadrant defined for diagrams $\mG: \FF\mD\to \Ab$ by the formula
\begin{equation}\label{spscog1}
 Ext^p_{\Ab^{\FF\mD}}(Lan^{\FF(f)}_q\DD\mC, \mG) \Rightarrow 
 Ext^n_{\Ab^{\FF\mC}}(\DD\mC, \mG\circ\FF(f))
\end{equation}
\end{proposition}
{\sc Proof.} It suffices to substitute into the spectral sequence of Proposition \ref{oberst1967} the diagram $\DD\mC$ instead of $\mF$, the functor $\FF(f)$ instead of $f$, $\FF\mC$ instead of $\mC$ and $\FF\mD$ instead of $\mD$.
\hfill$\Box$

Applying the Lemma \ref{degABN} to the spectral sequence of Proposition \ref{spscog} leads to the following assertion
\begin{proposition}\label{step1}
Let $N\geq 1$ be a natural number.
Edge homomorphisms
\begin{equation}\label{bouhom}
Ext^n(Lan^{\FF(f)}\DD\mC, \mG)\to Ext^n_{\Ab^{\FF\mC}}(\DD\mC, \mG\circ\FF(f))\end{equation}
spectral sequence (\ref{spscog1}) for all $\mG\in \Ab^{\FF\mD}$ are isomorphisms as $n\leq N$ and monomorphisms as $n=N+1$ if and only if $Lan^{\FF(f)}_q \DD\mC= 0$.
\end{proposition}

\begin{lemma}\label{step2}
The homomorphism $H^0_{\DF}(\mD, \mG)\to H^0_{\DF}(\mC, \mG\circ\FF(f))$ is 
invertible for every $\mG\in \Ab^{\FF\mD}$ if and only if the homomorphism $Lan^{\FF(f)}\DD\mC\to \DD\mD$ is an isomorphism.
\end{lemma}
{\sc Proof.} This homomorphism is an isomorphism if and only if
  $Hom(\DD\mD, \mG)\to Hom(\DD\mC, \mG\circ \FF(f))$ is an isomorphism.
Since 
$$
Hom(\DD\mC, \mG\circ \FF(f)) \cong Hom(Lan^{\FF(f)}\DD\mC, \mG),
$$
then this condition leads to the invertibility of the homomorphism 
$$
Hom(\DD\mD, \mG)\to Hom(Lan^{\FF(f)}\DD\mC, \mG)
$$ 
for all $\mG$, which is equivalent to the invertibility of the morphism $Lan^{\FF(f)}\DD\mC\to \DD\mD$.\hfill$\Box$

\begin{theorem}
Let $f: \mC\to \mD$ be a functor between small categories and $N\geq 1$ be a positive integer.
The following properties of the functor $f$ are equivalent
\begin{enumerate}
\item \label{prop1} For every natural system $\mG: \FF\mD\to \Ab$ and $n\geq 0$ the canonical homomorphisms of $H^n_{\DF}(\mD, \mG)\to H^n_{\DF}(\mC, \mG\circ\FF(f))$ are isomorphisms for all non-negative integers $n\leq N$ and monomorphisms for $n=N+1$.
\item \label{prop2} The homomorphism $Lan^{\FF(f)}\DD\mC\to \DD\mD$ is an isomorphism and for all $q$ in the interval $1\leq q\leq N$ the values of $Lan^{\FF(f)}_q\DD\mC$ are equal to $0$.
\end{enumerate}
\end{theorem}
{\sc Proof.} Consider any $N\geq 1$.
Since the property \ref{prop1} involves an isomorphism for $n=0$, and property 
\ref{prop2} is
isomorphism $Lan^{\FF(f)}\DD\mC\xrightarrow{\cong} \DD\mD$, then by Lemma \ref{step2} both isomorphisms can be included in each of the properties \ref{prop1} and \ref{prop2}.
By definition of cohomology, $H^n_{\DF}(\mC, \mG\circ \FF(f))=Ext^n(\DD\mC, \mG\circ\FF(f))$ and $H^n_{\DF}(\mD, \mG) = Ext^n(\DD\mD, \mG)$.
The isomorphism $Lan^{\FF(f)}\DD\mC\to \DD\mD$ gives isomorphisms
$$
H^n_{\DF}(\mD, \mG) = Ext^n(\DD\mD, \mG) \to
Ext^n(Lan^{\FF(f)}\DD\mC, \mG),
$$
  for all $n\geq 0$.
Proposition \ref{step1} implies the equivalence of the properties \ref{prop1} and \ref{prop2}.

\section{Conclusion}

An analysis of various theories of cohomology of small categories has shown that in most cases adjoint functors do not preserve the cohomology of a category when passing to the inverse image of a functor between small categories.
Nevertheless, it is concluded that the problem of characterizing functors between small cohomology-preserving categories has a solution for a large class of cohomology theories. This class includes the cohomologies constructed as derivatives of the limit functor, as well as the Baues-Wirsching, Hochschild-Mitchell, Thomason cohomology. This class consists of the compatible cohomologies introduced by us, defined on objects and morphisms belonging to Grothendieck constructions.

For each of these types of cohomology of small categories, we have obtained a means that allows us to determine whether the functor between small categories preserves the cohomology of the category when passing to the inverse image of this functor.
The verification of cohomology conservation conditions reduces to the calculation of integral homology groups of simplicial sets.


\begin{thebibliography}{1}

\bibitem{bau1985}
H.-J. Baues, G. Wirsching, 
\textit{Cohomology of small categories},
J. Pure Appl. Algebra {\bf 38} (1985), 187-211.

\bibitem{ber2021}
J. Bergfalk, 
\textit{The first omega alephs: From simplices to trees of trees to higher walks}, 
Adv. Math. {\bf 393} (2021), 108083.


\bibitem{bor1994}
F. Borceux, \textit{Handbook of Categorical Algebra 1. Basic Category Theory},
Cambridge: Cambridge University Press, 1994.


\bibitem{gal2012}
I. G\'alvez-Carrillo, F. Neumann, A. Tonks, \textit{Andr\'e spectral sequences
for Baues-Wirsching cohomology of categories}, J. Pure Appl. Algebra {\bf 216}
 (2012), 2549-2561.

\bibitem{gal2013}
I. G\'alvez-Carrillo, F. Neumann, A. Tonks, \textit{Thomason cohomology of
categories}, J. Pure Appl. Algebra {\bf 217} (2013), 2163-2179.

\bibitem{gal2021}
I. G\'alvez-Carrillo, F. Neumann, A. Tonks, \textit{Gabriel-Zisman cohomology 
and spectral sequences}, Appl. Categor. Struct. {\bf 29} (2021), 69-94.

\bibitem{gab1967}
P. Gabriel, M. Zisman, 
\textit{Calculus of fractions and homotopy theory}, 
Berlin-Heidelberg-New York: Springer-Verlag, 1967.

\bibitem{gro1957}
A. Grothendieck, \textit{Sur quelques points d'alg\`ebre homologique}. 
T\^ohoku Mathematical Journal, (2). {\bf 9}: 2 (1957), 119-221.

\bibitem{jib2006}
M. Jibladze, T. Pirashvili,  
\textit{Quillen cohomology and Baues-Wirsching cohomology of algebraic, theories}, Cahiers de topologie et g\'eom\'etrie 
diff\'erentielle cat\'egoriques, {\bf 47}:3 (2006), 163-205.

\bibitem{X1997}
A. A. Khusainov, {\textit Comparing dimensions of a small category},
Sib. Math. J. {\bf 38} (1997), 1230-1240.
    
\bibitem{mac2004}
S. Mac Lane, 
\textit{Categories for the Working Mathematician}. 
New York: Springer-Verlag (1998).

\bibitem{mac1963}
S. Mac\,Lane, 
{\it Homology},
Die Grundlehren der mathematischen Wissenschaften {\bf 114},
Berlin-Heidelberg-New York: Springer-Verlag (1975).

\bibitem{mit1972}
B. Mitchell, Rings with several objects, Adv. Math. {\bf 8} (1972) 1-161.

\bibitem{mit1973}
B. Mitchell,
     \textit{The cohomological dimension of a directed set}, 
Canad. J. Math. {\bf 25}: 2 (1973), 233-238.

\bibitem{mur2006}
F. Muro, \textit{On the functoriality of cohomology of categories}, 
J. Pure Appl. Algebra {\bf 204} (2006) 455 - 472.

\bibitem{obe1967}
U. Oberst,
     \textit{Basisweiterung in der Homologie kleiner
     Kategorien}. Math. Z. {\bf 100} (1967), 36-58.

\bibitem{obe1968}
U. Oberst,
\textit{Homology of categories and exactness
       of direct limits}, Math. Z.  {\bf 107} (1968), 87-115.
  
\bibitem{oli1994}       
B. Oliver,
\textit{Higher limits via Stenberg representation},
Communication in Algebra. {\bf 22}:4 (1994), 1381-1393.

\bibitem{pir2006}
T. Pirashvili, M.J. Redondo, Cohomology of the Grothendieck construction, Manuscr. Math. {\bf 120} (2006), 151-162.

\bibitem{qui1973}
D. G. Quillen,
 \textit{Higher algebraic K-theory: I}, Higher K-Theories.
 Lecture Notes in Math., {\bf 341}.
Berlin-Heidelberg: Springer-Verlag (1973), 85-147.

\bibitem{tho1979}
R. W. Thomason,
{\it Homotopy colimits in the category of small categories.} 
Math. Proc. Camb. Phil. Soc.  {\bf 85}:1 (1979), 91-109.

\bibitem{yal2023}
E. Yal\c{c}{\i}n, 
\textit{LHS-spectral sequences for regular extensions of categories.}, 
New York, 2023. 43 p.
(Preprint, Cornell Univ.); https://arxiv.org/abs/2305.02000

\bibitem{yal2022}
E. Yal\c{c}{\i}n, 
\textit{Higher limits over the fusion orbit category}, 
Adv. Math. {\bf 406} (2022), 108482.


\bibitem{nlabGC}
https://ncatlab.org/nlab/show/Grothendieck+construction\#Definition



\end{thebibliography}
\end{document}